\documentclass[11pt]{amsart}
\usepackage[hmargin=2.5cm,vmargin=2.5cm]{geometry}
\usepackage{amsfonts}
\usepackage{amsthm}
\usepackage[english]{algorithm2e}
\usepackage{hyperref}
\usepackage{amsmath}
\usepackage{amssymb}
\usepackage[T1]{fontenc}
\newtheorem{Theorem}{Theorem}[section]
\newtheorem{Definition}[Theorem]{Definition}
\newtheorem{Remark}[Theorem]{Remark}
\newtheorem{Example}[Theorem]{Example}

\newtheorem{Proposition}[Theorem]{Proposition}
\newtheorem{Corollary}[Theorem]{Corollary}

\title{Hom-Lie-Yamaguti Superalgebras}
\author[Sylvain Attan and Donatien Gaparayi.]
       { Donatien Gaparayi $^1$ and Sylvain Attan $^2$
       \\\\
          $^{1}$  Ecole Normale Sup\'erieure (E.N.S),\\
        BP 6983 Bujumbura, Burundi\\
       $^2$ D\'epartement de Math\'ematiques\\
       Universit\'e d'Abomey-calavi\\
       01 BP 4521, Cotonou 01, B\'enin\\
       $^1$gapadona@yahoo.fr, $^2$ syltane2010@yahoo.fr
       }
\begin{document}
\maketitle
\begin{abstract}
(Multiplicative) Hom-Lie-Yamaguti superalgebras which generalize Hom-Lie supertriple  systems (and
 subsequently ternary multiplicative Hom-Nambu superalgebras) and  Hom-Lie superalgebras in the same way as Lie-Yamaguti superalgebras
 \cite{Frac} generalize Lie supertriple systems and Lie superalgebras are defined. We show  that the category of (multiplicative)  
 Hom-Lie-Yamaguti superalgebras is closed under twisting by self-morphisms. Construction of some examples of Hom-Lie-Yamaguti superalgebras
 is given. The notion of an $n^{th}-$derived (binary) Hom-superalgebras is extended to the one of an $n^{th}-$derived binary-ternary 
 Hom-superalgebras and it is shown that the category of Hom-Lie-Yamaguti superalgebras is closed under the process of taking $n^{th}-$derived 
 Hom-superalgebras. 
\end{abstract}

{\bf Keywords:} Lie-Yamaguti superalgebra \cite{Frac} (i.e.generalized Lie supertriple system, Lie superalgebra), Hom-Lie-Yamaguti
superalgebra (i.e.generalized Hom-Lie supertriple system, Hom-Lie superalgebra).
\section{Introduction}
\cite{Frac} A Lie-Yamaguti superalgebra is a triple $(L,\ast,\{,,\})$ where 
 $L = L_0 \oplus L_1$ is $\mathbb{K}-$vector superspace i.e
 $\mathbb{Z}_2-$graded vector space, $"\ast"$ a bilinear map (the binary superoperation on $L$) and $"\{,,\}"$ a trilinear map (the ternary superoperation 
 on $L$)  such that
\begin{eqnarray*}
&& (SLY1) \,\,\,x\ast y = - (-1)^{|x||y|}y\ast x,\\
&& (SLY2)\,\, \{x,y,z\} = -(-1)^{|x||y|}\{y,x,z\},\\ 
&& (SLY3)\,\,\{x,y,z\} + (-1)^{|x|(|y|+|z|)}\{y,z,x\} + (-1)^{|z|(|x|+|y|)}\{z,x,y\} + J(x,y,z)= 0,\\
&& (SLY4)\,\, \{x\ast y,z,u\} + (-1)^{|x|(|y|+|z|)}\{y\ast z,x,u\} + (-1)^{|z|(|x|+|y|)}\{z\ast x,y,u\} = 0,\\
&& (SLY5)\,\, \{x,y,u\ast v\} = \{x,y,u\}\ast v + (-1)^{|u|(|x|+|y|)}u\ast \{x,y,v\},\\
&& (SLY6)\,\, \{x,y,\{u,v,w\}\} = \{\{x,y,u\},v,w\} + (-1)^{|u|(|x|+|y|)}\{u,\{x,y,v\},w\} \\
&& \hspace{5cm} + (-1)^{(|x|+|y|)(|u|+|v|)}\{u,v,\{x,y,w\}\},
\end{eqnarray*}
for all $x,y,z,u,v,w$ in $L$ and 
\begin{equation}\label{SJ}
J(x,y,z) = (x\ast y)\ast z +(-1)^{|x|(|y|+|z|)}(y\ast z)\ast x + (-1)^{|z|(|x|+|y|)}(z\ast x)\ast y.
\end{equation}
The relation (\ref{SJ}) is called super-Jacobian.
Observe that if $x\ast y = 0$, for all $x,y$
in $L$, then a Lie-Yamaguti superalgebra $(L,\ast, \{, , \})$ reduces to a Lie supertriple system $(L, \{, , \})$ as defined in \cite{Okubo2} 
and if $\{x,y,z\} = 0$ for all $x, y, z$ in $L$, then $(L, \ast, \{, , \})$ is a Lie superalgebra $(L,\ast)$ \cite{Okubo2}. 
Recall that a Lie supertriple system \cite{Okubo2} is a pair $(T,\{,,\})$ where $T = T_0 \oplus T_1$ is a  
$\mathbb{K}-$vector superspace  and $"\{,,\}"$ a trilinear map (the ternary superoperation on $T$) such that
\begin{eqnarray*}
&& (i)\,\,\{x,y,z\} = (-1)^{|x||y|}\{y,x,z\}\\
&& (ii)\,\,\,\{x,y,z\} + (-1)^{|x|(|y|+|z|)}\{y,z,x\} + (-1)^{|z|(|x|+|y|)}\{z,x,y\} = 0\\
&& (iii)\,\,\, \{x,y,\{u,v,w\}\} = \{\{x,y,u\},v,v\} + (-1)^{(|x|+|y|)|u|}\{u,\{x,y,v\},w\}\\
&& \hspace{4.3cm} + (-1)^{(|x|+|y|)(|u|+|v|)}\{u,v,\{x,y,w\}\}
\end{eqnarray*}
for all $x, y, z$ in $T$ and a Lie superalgebra \cite{Okubo2} is a pair $(A,\ast)$ where $A= A_0 \oplus A_1$ is a 
$\mathbb{K}-$vector superspace with $"\ast"$ a bilinear map (the binary superoperation on $A$) such that
\begin{eqnarray*}
&& (i) \,\,\,x\ast y = (-1)^{|x||y|}y\ast x\\
&& (ii)\,\,\, J(x,y,z) = 0
\end{eqnarray*}
 for all $x, y, z$ in $A$.\newline

A Hom-type generalization of a kind of algebra is obtained by certain twisting of the defining identies by a linear self-map, called
the twisting map, in such way that when the twisting map is identity map, then one recovers the original kind of algebra. In this scheme,
e.g., associative algebras and Leibniz algebras are twisted into Hom-associative algebras and Hom-Leibniz algebras respectively \cite{MAK3}
and, likewise, Hom-type analogues of Novikov algebras, alternative algebras, Jordan algebras or Malcev algebras are defined and discussed
in \cite{MAK1},\cite{YAU3}, \cite{YAU4}. Also, Hom-Lie triple systems \cite{YAU5}, $n-$ary Hom-Nambu and Hom-Nambu-Lie algebras  
\cite{Ataguema1}, were introduced in. In the some way, this generalization of binary or ternary algebras has been extended on the 
binary-ternary algebras. Indeed, Hom-Akivis algebras \cite{NOU1}, Hom-Lie-Yamaguti algebras \cite{GAP} and Hom-Bol algebras \cite{ATTAN} 
which generalize Akivis algebras, Lie-Yamaguti algebras and Bol algebras respectively are introduced. One could say that the theory of 
Hom-algebras originated in \cite{HAR1} (see also \cite{HAR2},\cite{HAR3}) in the study of deformations of Witt and Virasoro algebras
(in fact, some $q-$deformations of Witt and Virasoro algebras have a structures of a Hom-Lie algebra \cite{HAR1}). Some algebraic abstractions
of this study are given in \cite{MAK3}, \cite{YAU2}. For further more information on other Hom-type algebras, one may refer to, e.g., 
\cite{Ataguema1}, \cite{Gohr2}, \cite{MAK1}, \cite{YAU3}, \cite{YAU4}, \cite{YAU5}.\newline 

In \cite{MAK5}, the authors introduce Hom-associative superalgebras and Hom-Lie superalgebras which generalize associative superalgebras
\cite{Okubo2} and Lie superalgebras \cite{Okubo2} respectively. They provide a  way for constructing Hom-Lie superalgebras from
Hom-associative superalgebras which extend the fundamental construction of Lie superalgebras from associative superalgebras via 
supercommutator bracket  (see the Proposition 1.1. in \cite{Okubo2}). Indeed, they show also that the supercommutator bracket defined using 
the multiplication in a Hom-associative superalgebra leads naturally to Hom-Lie superalgebras. 
In \cite{MAK4} the authors introduce the Hom-alternative, Hom-Malcev and Hom-Jordan superalgebras wich are the generalization of Alternative,
Malcev and Jordan superalgebras respectively. \newline

Our present study extends the Hom-type generalization of binary superalgebras to the one of ternary superalgebras or binary-ternary 
superalgebras. The purpose of this paper is to introduced  Hom-type generalization of Lie-Yamaguti superalgebras \cite{Frac}, called 
Hom-Lie-Yamaguti superalgebras. It is also to extend the notion of an $n^{th}-$derived (binary) Hom-superalgebras \cite{MAK4} to the one
of an $n^{th}-$derived ternary or binary-ternary Hom-superalgebras and  we shown that ternary or binary-ternary superalgebras are closed 
under the process of taking $n^{th}-$derived Hom-superalgebras.\newline
 
 The rest of this paper is organized as follows. In Section two, we recall basic definitions in Hom-superalgebras theory and useful 
 results  about Hom-associative superalgebras and Hom-Lie superalgebras. In \cite{MAK5}, the autors show that the supercommutator bracket
 defined using the multiplication in a Hom-associative superalgebra leads naturally to Hom-Lie superalgebra. Also, we recall the notion
 of $n^{th}-$derived (binary) Hom-superalgebra introduced in \cite{MAK5} and as example, we show that Hom-Lie superalgebras are closed 
 under the process of taking $n^{th}-$derived (binary) Hom-superalgebras (see the Proposition \ref{SH-Lie}). In the third Section, we 
 introduce ternary and binary-ternary  Hom-superalgebras. 
 In particulary, Hom-Lie-Yamaguti superalgebras which are binary-ternary Hom-superalgebras are defined. 
 Hom-Lie-Yamaguti superalgebras  generalize Hom-Lie triple  supersystems (and  subsequently ternary multiplicative Hom-Nambu superalgebras)
 and  Hom-Lie superalgebras. We provide that any non-Hom-superassociative Hom-superalgebra is a Hom-supertriple system. We also extend a
 Yau's theorem  \cite{YAU2}  on Hom-Lie algebras to Hom-Lie-Yamaguti superalgebras and we conclude that these category of
 Hom-superalgebras are closed   under twisting by self-morphisms (see Theorem \ref{Thm SHLY}). Some examples of
 Hom-Lie-Yamaguti superalgebras are constructed using Theorem \ref{Thm SHLY} via Corollary \ref{Metho}. In the last  Section, we extend the
 notion of  $n^{th}-$ derived (binary) Hom-superalgebra to the case of ternary and binary-ternary Hom-superalgebras.  We show that the category  of  
 Hom-Lie-Yamaguti superalgebras are closed under the process of taking $n^{th}-$derived Hom-superalgebras.
\section{Some basics on superalgebras}
We recall some basic facts about Hom-superalgebras, including Hom-associative and Hom-Lie superalgebras \cite{MAK5}. Also, we recall the
 notion of  (binary) $n^{th}-$derived Hom-superalgebras. As example, we show that the category of Hom-Lie superalgebras is closed under the
 process of  taking $n^{th}-$derived Hom-superalgebras \cite{MAK5}. \newline
 
\begin{Definition}
\begin{itemize}
\item  [$(i)$] Let $f:(A, \ast, \alpha) \rightarrow (A' , \ast',\alpha')$ be a map, where $A = A_0\oplus A_1$ and $A' = A'_0\oplus A'_1$ are 
$\mathbf{Z}_2-$graded vector spaces. The map $f$ is called an even (resp. odd ) map if  $f(A_i) \subset A'_i$ (resp. $f (A_i ) \subset A'_{i+1})$,
for $i = 0,1.$
\item [$(ii)$] A \textbf{Hom-superalgebra} is a triple $(A,\ast, \alpha)$ in which $A = A_0\oplus A_1$ is a $\Bbb K-$super-module, 
$\ast:A \times A \rightarrow A$ is an even  bilinear map, and $\alpha: A \rightarrow A$ is an even linear map such that 
$\alpha(x \ast y) = \alpha(x) \ast \alpha(y)$ (multiplicativity).
\end{itemize}
% \newline
\end{Definition}

\begin{Remark}
For convenience, we assume throughout this paper that all Hom-superalgebras are multiplicative. 
\end{Remark}

\begin{Definition}
  Let $(A,\ast, \alpha)$ be a Hom-superalgebra, that is a $\Bbb K-$vector superspace $A$ together with a  multiplication $"\ast"$  and an
  even   linear  self-map $\alpha$. 
  \begin{itemize}
 \item [$(i)$]  \cite{MAK4} The Hom-superassociator of $A$ is the trilinear map $as_{\alpha}: A \times A \times A \rightarrow A$ defined as
\begin{equation}\label{asA}
as_{\alpha} = \ast \circ (\ast \otimes \alpha - \alpha \otimes \ast).
\end{equation}
In terms of elements, the map $as_{\alpha}$ is given by
\begin{equation*}
as_{\alpha}(x,y,z) = (x\ast y)\ast \alpha(z)- \alpha(x)\ast(y\ast z).
\end{equation*}
\item [$(ii)$] \cite{MAK5} The Hom-superalgebra $(A,\ast,\alpha)$ is called a Hom-associative superalgebra if $$as_{\alpha}(x,y,z) = 0,\,\,\, 
\forall\,\,x,y,z\in A$$
\item [$(iii)$] \cite{MAK4} The Hom-super-Jacobian of $A$ is the trilinear map $J_{\alpha} : A\times A\times A \rightarrow A$ defined
as $J_\alpha(x,y,z):= (x\ast y)\ast \alpha(z) + (-1)^{|x|(|y|+|z|)}(y\ast z)\ast \alpha(x) + (-1)^{|z|(|x|+|y|)}(z\ast x)\ast \alpha(y)$
for all $ x,y,z \in A.$ If $\alpha = Id$, then a Hom-super-Jacobian reduces to (\ref{SJ}).
\item [$(iv)$] \cite{MAK5} The Hom-superalgebra $(A,\ast,\alpha)$ is called a Hom-Lie superalgebra if
\begin{eqnarray}\label{HSJ}\nonumber
&& x\ast y = - (-1)^{|x||y|}y\ast x\\
&& J_{\alpha} (x,y,z) = 0 \,\,\,\mbox{(the Hom-super-Jacobi identity)}
\end{eqnarray}
for all $x,y,z\in A.$ 
If $\alpha = Id,$ a Hom-Lie superalgebra reduces to a usual Lie superalgebra.
\end{itemize}
\end{Definition}
Recall that the supercommutator bracket defined using the multiplication in any Hom-associative superalgebra leads naturally to a Hom-Lie 
superalgebra (see the Proposition 2.6 in \cite{MAK5}). Others Hom-Lie superalgebras can be constructed from Lie superalgebras using the 
Theorem 2.7 in \cite{MAK5}. In the following, we recall the notion of (binary) $n^{th}-$derived Hom-superalgebra.

\begin{Definition}\label{H−B−D}
\cite{MAK4} Let $(A,\ast,\alpha)$ be a Hom-superalgebra and $n\geq 0$ an integer ( $"\ast"$ is the binary operation on $A$). The
Hom-superalgebra $A^n$ defined by 
$$A^n := (A, \ast^{(n)},\alpha^{2^n}),\,\,\mbox{where} \,\,(x\ast^{(n)} y):= \alpha^{2^{n}-1}(x\ast y), \forall\,\, x, y \in A $$
is called the {\it $n^{th}-$derived Hom-superalgebra} of $A$.
\end{Definition}
For simplicity of exposition, $\ast^{(n)}$ is written as $ \ast^{(n)} = \alpha^{2^{n}-1}\circ \ast.$ Then notes that
$A^0 = (A,\ast,\alpha), A^1 = (A, \ast^{(1)} = \alpha \circ \ast,\alpha^2),$ and $A^{n+1} = (A^n)^1 $.\newline

\begin{Proposition}\label{SH-Lie}
 Let $(A,[,],\alpha)$ be a multiplicative Hom-Lie superalgebra. Then the $n^{th}-$derived Hom-superalgebra 
 $A^n = (A,[,]^{(n)} = \alpha^{2^n-1} \circ [,], \alpha^{2^n})$ is also a multiplicative Hom-Lie superalgebra for each $n\geq 0.$
\end{Proposition}

\textbf{Proof:} Observe that $[x,y]^{(n)} = - (-1)^{|x||y|}[y,x]^{(n)}.$ Indeed, 
$$[x,y]^{(n)} = \alpha^{2^{n}-1}([x,y]) = \alpha^{2^{n}-1}((-1)^{|x||y|}[y,x]) = (-1)^{|x||y|}\alpha^{2^{n}-1}([y,x]) = (-1)^{|x||y|}[y,x]^{(n)}.$$ 
Then, the superantisymmetry of $[x,y]^{(n)}$ holds in $A^n$. Next, we have
\begin{eqnarray*}
&& [[x,y]^{(n)},\alpha^{2^n}(z)]^{(n)} + (-1)^{|x|(|y|+|z|)}[[y,z]^{(n)},\alpha^{2^n}(x)]^{(n)} 
+ (-1)^{|x|(|y|+|z|)}[[z,x]^{(n)},\alpha^{2^n}(y)]^{(n)}\\
& = & \alpha^{2^{n}-1}([\alpha^{2^{n}-1}([x,y]),\alpha^{2^n}(z)]) + (-1)^{|x|(|y|+|z|)}\alpha^{2^{n}-1}([\alpha^{2^{n}-1}([y,z]),\alpha^2(x)])\\
& + & (-1)^{|z|(|x|+|y|)}\alpha^{2^{n}-1}([\alpha^{2^{n}-1}([z,x]),\alpha^{2^n}(y)])\\
& = & \alpha^{2(2^{n}-1)}([[x,y],\alpha^2(z)] + (-1)^{|x|(|y|+|z|)}[[y,z],\alpha^2(x)]+ (-1)^{|z|(|x|+|y|)}[[z,x],\alpha^2(y)])\\
& = & \alpha^{2(2^{n}-1)}(0)\quad \mbox{(by (\ref{HSJ}))}\\ 
& = & 0
\end{eqnarray*}
and so the Hom- superJacobi identity (\ref{HSJ}) holds in $A^{(n)}.$ Thus, we conclude that $A^{(n)}$ is a (multiplicative) Hom-Lie 
superalgebra. In section 4, we extend this notion of $n^{th}-$one derived (binary) Hom-superalgebra to the case of ternary and 
binary-ternary Hom-superalgebras.%\newline
%%%%%%%%%%%%%%%%%%%%%%%%%%%%%%%%%%%%%%%%%%%%%%%%%%%%%%%%%%%%%%%%%%%%%%%%%%%%%%%%%%%%%%%%%%%%%%%%%%%%%%%%%%%%%%%%%%%%%%%%%%%%%%%%%%%%%%%%%%
\section{Ternary and binary-ternary Hom-superalgebras}
%%%%%%%%%%%%%%%%%%%%%%%%%%%%%%%%%%%%%%%%%%%%%%%%%%%%%%%%%%%%%%%%%%%%%%%%%%%%%%%%%%%%%%%%%%%%%%%%%%%%%%%%%%%%%%%%%%%%%%%%%%%%%%%%%%%%%%%%%%%%
In this section, we introduce Hom-Lie supertriple  systems and consequently ternary multiplicative Hom-Nambu superalgebras.  
Hom-Lie-Yamaguti superalgebras which are binary-ternary Hom-superalgebras are defined. Hom-Lie-Yamaguti superalgebras reduce to Lie-Yamaguti 
superalgebras \cite{Frac} when the twisting map is identity map. Note also that Hom-Lie-Yamaguti superalgebras generalize Hom-Lie supertriple
systems (and  subsequently ternary multiplicative Hom-Nambu superalgebras) and  Hom-Lie superalgebras. We finish this section by 
giving some examples of Hom-Lie-Yamaguti superalgebras constructed by using Theorem \ref{Thm SHLY} via Corollary \ref{Metho}.%\newline
%%%%%%%%%%%%%%%%%%%%%%%%%%%%%%%%%%%%%%%%%%%%%%%%%%%%%%%%%%%%%%%%%%%%%%%%%%%%%%%%%%%%%%%%%%%%%%%%%%%%%%%%%%%%%%%%%%%%%%%%%%%%%%%%%%%%%%%%%%
\subsection{Definitions and Theorem of construction}
%%%%%%%%%%%%%%%%%%%%%%%%%%%%%%%%%%%%%%%%%%%%%%%%%%%%%%%%%%%%%%%%%%%%%%%%%%%%%%%%%%%%%%%%%%%%%%%%%%%%%%%%%%%%%%%%%%%%%%%%%%%%%%%%%%%%%%%%%%
In this subsection, we give some definitions of ternary Hom-superalgebras and we show that any non-Hom-superassociative
Hom-superalgebra has a natural structure of Hom-supertriple system. In the Theorem \ref{Thm SHLY}, we show 
that Hom-Lie-Yamaguti superalgebras are closed under twisting by self-morphisms respectively.
\begin{Definition}
  A \textbf{ternary Hom-superalgebra} is triple $(T,\{,,\},\alpha = (\alpha_1,\alpha_2))$ constitued by a vector $\Bbb K-$superspace
  $"T = T_0 \oplus T_1"$, a   trilinear map   $\{,,\} : T\times T\times T \longrightarrow T,$ and even linear maps $\alpha_i : T\rightarrow T,
  \,\,i = 1,2,$   called the \textbf{twisting   maps}. The algebra $(T,\{,,\},\alpha = (\alpha_1,\alpha_2))$ is said \textbf{multiplicative} 
  if   $\alpha_1 = \alpha_2 := \alpha$ and $\alpha (\{x,y,z\}) = \{\alpha (x),\alpha (y),\alpha (z)\}$ for all $x, y, z \in T $.\newline
  \end{Definition}

\begin{Definition}
 A ({\it multiplicative}) ternary Hom-Nambu superalgebra is (a multiplicative) ternary Hom-superalgebra $(T,\{,,\},\alpha)$ satisfying
\begin{eqnarray}\label{Nambu1} \nonumber
\{\alpha (x), \alpha (y), \{u,v,w\}\}
& = & \{\{x,y,u\},\alpha (v), \alpha(w)\}\\\nonumber
& + & (-1)^{|u|(|x|+|y|)}\{\alpha (u), \{x,y,v\}, \alpha (w)\}\\
& + & (-1)^{(|x|+|y|)(|u|+|v|)}\{\alpha (u),\alpha (v),\{x,y,w\}\}
\end{eqnarray}
%}
 for all $u,v,w,x,y,z \in T.$ The condition (\ref{Nambu1}) is called the {\it ternary Hom-Nambu superidentity.} In general, the ternary
 Hom-Nambu superidentity reads:
\begin{eqnarray*}\label{Nambu2}
\{\alpha_1(x), \alpha_2(y), \{u,v,w\}
& = & \{\{x,y,u\},\alpha_1(v), \alpha_2(w)\}\\
& + & (-1)^{(|x|+|y|)|u|}\{\alpha_1(u), \{x,y,v\}, \alpha_2(w)\}\\ 
& + & (-1)^{(|x|+|y|)(|u|+|v|)}\{\alpha_1(u),\alpha_2(v),\{x,y,w\}\}
\end{eqnarray*}
for all $u,v,w,x,y,z \in T,$ where $\alpha_1$ and $\alpha_2$ are linear  self-maps of $T.$ \newline
\end{Definition}

\begin{Definition}\label{H-STS}
 A (multiplicative) Hom-supertriple system is a (multiplicative) ternary Hom-superalgebra  $(T,\{.,.,.\},\alpha)$ such that 
 \begin{itemize}
  \item [$(i)$] $\{x,y,z\} = - (-1)^{|x||y|}\{y,x,z\},$
  \item [$(ii)$] $\{x,y,z\} + (-1)^{|x|(|y|+|z|)}\{y,z,x\} + (-1)^{|z|(|x|+|y|)}\{z,x,y\} = 0\,\,\, \forall\,\,x,y,z \in T.$
 \end{itemize}
\end{Definition}
\begin{Remark} 
 Our definition here is motivated by the concern of giving a Hom-type analogue of the relationships between nonassociative superalgebras and 
 supertriple systems (see Proposition \ref{Prop-STS}). \newline
\end{Remark}

 \begin{Proposition}\label{Prop-STS}
  Any non-Hom-associative superalgebra is a Hom-supertriple system.
 \end{Proposition}

 \textbf{Proof}: Let $(A,\ast,\alpha)$ be a non-Hom-associative superalgebra. Define the supercommutator  by 
 $[x,y] := x \ast y - (-1)^{|x||y|} y \ast x$  and the ternary superoperation by 
 $\{x,y,z\} := [[x, y], \alpha(z)] - as_{\alpha}(x, y, z) + (-1)^{|x||y|}as_{\alpha}(y, x, z)$
 for all homogeneous elements $x,y,z \in A$ and where the superassociator $as(x,y,z)$ is given by (\ref{asA}).
 Then we have $\{x,y,z\} = - (-1)^{|x||y|}\{y,x,z\}$ and  $\circlearrowleft_{(x,y,z)} (-1)^{|x||z|}\{x,y,z\} = 0.$ Thus  $(A, \{,,\},\alpha)$
 is a Hom-supertriple system.   \hfill$\square$% \newline

\begin{Remark}
 For $\alpha = Id$ (the identity map), we recover the supertriple system with ternary superoperation 
 $\{x,y,z\} = [[x,y],z]-as(x,y,z) + (-1)^{|x||y|}as(y,x,z)$ that is associated to each nonassociative superalgebra and where $[x,y]$ and
 $as(x,y,z)$ are supercommutator and superassociator respectively for all homogeneous elements $x,y,z.$% \newline
\end{Remark}

\begin{Definition}\label{LTSup}
 A \textbf{Hom-Lie supertriple system} {\it is a Hom-supertriple system $(A,\{,,\},\alpha )$ satisfying  the ternary Hom-Nambu superidentity 
 (\ref{Nambu1}).}  When $\alpha = Id$, a Hom-Lie supertriple system reduces to a Lie supertriple system.%\newline
\end{Definition}
One can note that  Hom-Bol superalgebras introduced in \cite{ANI} may be viewed as some generalization of  Hom-supertriple systems.

In the following, we now give the definition of the basic object of this paper and we consider construction methods for Hom-LY
superalgebras. These methods allow to find examples of Hom-LY superalgebras starting from ordinary LY superalgebras or even from Malcev 
superalgebras. Recall that from a Malcev superalgebra $(A,\ast)$ and consider on $A$ the ternary superoperation 
\begin{equation}\label{M-LY}
\{x,y,z\} := x \ast (y\ast z) - (-1)^{|x||y|}y\ast (x\ast z) + (x\ast y)\ast z, \forall\,\,x,y,z \in A,
\end{equation}
then $(A,\ast, \{,,\})$ is a Lie-Yamaguti superalgebra \cite{Frac}. First, as the main tool, we show that the category of (multiplicative)
Hom-LY superalgebras is closed under self-morphisms (see the Theorem \ref{Thm SHLY}). %\newline

\begin{Definition}\label{SHLY-Def}
 {\it A \textbf{Hom-Lie Yamaguti superalgebra} (Hom-LY superalgebra for short) is a quadruple $(L,\ast,\{,,\},\alpha )$ 
in which $L$ is $\Bbb K$-vector superspace, $"\ast"$ a binary superoperation and $"\{,,\}"$ a ternary superoperation on $L$, and 
$\alpha : L\rightarrow L$ an even linear map such that %\newline
\begin{eqnarray*}
&& (SHLY1) \,\,\alpha (x \ast y) = \alpha(x)\ast \alpha(y),\\
&& (SHLY2) \,\,\alpha (\{x,y,z\}) = \{\alpha(x,)\alpha(y),\alpha(z)\},\\ 
&& (SHLY3)\,\, x\ast y = - (-1)^{|x||y|}y\ast x,\\ 
&& (SHLY4) \,\,\{x,y,z\} = - (-1)^{|x||y|}\{y,x,z\},\\
&& (SHLY5) \,\,\circlearrowleft_{(x,y,z)} (-1)^{|x||z|}[(x\ast y)\ast \alpha(z) + \{x,y,z\}] = 0,\\ 
&& (SHLY6)\,\, \circlearrowleft_{(x,y,z)}(-1)^{|x||z|}[\{x\ast y, \alpha(z),\alpha(u)\}] = 0,\\
&& (SHLY7) \,\,\{\alpha(x),\alpha(y),u\ast v\} = \{x,y,u\} \ast \alpha^2(v) + (-1)^{|u|(|x|+|y|)}\alpha^2(u)\ast \{x,y,v\},\\ 
&& (SHLY8)\,\, \{\alpha^2(x),\alpha^2(y),\{u,v,w\}\} =  \{\{x,y,u\},\alpha^2(v),\alpha^2(w)\}\\
 && \hspace{6.5cm}+\,\,  (-1)^{|u|(|x|+|y|)}\{\alpha^2(u),\{x,y,v\},\alpha^2(w)\}  \\
 && \hspace{6.5cm}+ \,\, (-1)^{(|x|+|y|)(|u|+|v|)}\{\alpha^2(u),\alpha^2(v),\{x,y,w\}\},
\end{eqnarray*}
for all $u,v,w,x,y,z \in L$ and where $\circlearrowleft_{(x,y,z)}$ denotes the sum over cyclic permutation of $x,y,z.$ } \newline
\end{Definition}
Note that the conditions {$(SHLY1)$} and {$(SHLY2)$} mean the multiplicativity of $(L,\ast,\{,,\},\alpha ).$\newline

\begin{Remark} 
 
\begin{itemize} 
\item [$(1)$] If $\alpha = Id,$ then the Hom-LY superalgebra $(L,\ast,\{,,\},\alpha )$ reduces to a LY superalgebra $(L,\ast,\{,,\})$ 
(see $(SLY1)-(SLY6)$). 
\item [$(2)$] If $x\ast y = 0,$ for all $x,y \in L,$ then $(L,\ast,\{,,\},\alpha )$ becomes a Hom-Lie supertriple system  $(L,\{,,\},\alpha^2 )$
and, subsequently, a ternary Hom-Nambu superalgebra (since, by Definition \ref{LTSup}, any Hom-Lie supertriple system is automatically a 
ternary Hom-Nambu superalgebra). 
\item [(3)] If $\{x,y,z\} = 0$ for all $x,y,z\in L,$ then the Hom-LY superalgebra $(L,\ast,\{,,\},\alpha )$ becomes a Hom-Lie superalgebra 
$(L,\ast,\alpha ).$%\newline
\end{itemize}
\end{Remark} 

 \begin{Theorem}\label{Thm SHLY}
 Let $A_\alpha := (A,\ast,\{,,\},\alpha )$ be a Hom-LY superalgebra and let $"\beta"$ be an even endomorphism  of the superalgebra 
 $(A,\ast,\{,,\})$  such that $\beta\alpha = \alpha\beta$.  Let $\beta^0 = Id$ and, for any $n\ge1$ $\beta^n := \beta \circ \beta^{n-1}$.
 Define on $A$
 the superoperations \\ $x\ast_\beta y := \beta^n(x\ast y),$\\  $\{x,y,z\}_\beta := \beta^{2n}(\{x,y,z\}),$ \\for all $x,y,z \in A$.
 Then, $A_{\beta^n} := (A,\ast_\beta,\{,,\}_\beta,\beta^n\alpha )$ are Hom-LY superalgebras, with $n\ge1$.
 \end{Theorem}

 {\bf Proof.} First, we observe that the condition $\beta\alpha = \alpha\beta$ implies $\beta^n\alpha = \alpha\beta^n, n\ge1.$ Next we have \\
 $(\beta^n\alpha)(x\ast_\beta y)$ =  $(\beta^n\alpha)(\beta^n(x)\ast \beta^n(y))$
 = $\beta^n((\alpha\beta^n)(x)\ast (\alpha\beta^n)(y))$\\
 = $(\alpha\beta^n)(x)\ast_\beta (\alpha\beta^n)(y)$
 = $(\beta^n\alpha)(x)\ast_\beta (\beta^n\alpha)(y)$ and we get $(SHLY1)$ for $A_{\beta^n}$. Likewise, the condition $\beta\alpha = \alpha\beta$ 
   implies $(SHLY2)$. The identities $(SHLY3)$ and $(SHLY4)$ for $A_{\beta^n}$ follow from  the skew-supersymmetry of $"\ast"$ and $"\{,,\}"$ 
   respectively.

  Consider  now  $\circlearrowleft_{x,y,z}[(x\ast_\beta y)\ast_\beta (\beta^n\alpha)(z) + \{x,y,z\}_\beta].$  Then
\begin{eqnarray*}
&&  (x\ast_\beta y)\ast_\beta (\beta^n\alpha)(z) + (-1)^{|x|(|y|+|z|)}(y\ast_\beta z)\ast_\beta (\beta^n\alpha)(x)
  +  (-1)^{|z|(|x|+|y|)}(z\ast_\beta x)\ast_\beta (\beta^n\alpha)(y)\\\nonumber
&& + \{x,y,z\}_\beta + (-1)^{|x|(|y|+|z|)}\{y,z,x\}_\beta + (-1)^{|z|(|x|+|y|)}\{z,x,y\}_\beta \\\nonumber
&& =  \beta^n(\beta^n(x\ast y)\ast (\beta^n(\alpha(z))) +(-1)^{|x|(|y|+|z|)}\beta^n(\beta^n(y\ast z)\ast (\beta^n(\alpha(x)))\\\nonumber
&& + (-1)^{|z|(|x|+|y|)}\beta^n(\beta^n(z\ast x)\ast (\beta^n(\alpha(y))) +  \beta^{2n}(\{x,y,z\}) + (-1)^{|x|(|y|+|z|)}\beta^{2n}(\{y,z,x\}\\\nonumber
&& + (-1)^{|z|(|x|+|y|)}\beta^{2n}(\{z,x,y\} \\\nonumber
&& =  \beta^{2n}((x\ast y)\ast \alpha(z)) +(-1)^{|x|(|y|+|z|)}\beta^{2n}((y\ast z)\ast \alpha(x))\\\nonumber
&& + (-1)^{|z|(|x|+|y|)}\beta^{2n}((z\ast x)\ast \alpha(y)) +  \beta^{2n}(\{x,y,z\}) + (-1)^{|x|(|y|+|z|)}\beta^{2n}(\{y,z,x\}\\\nonumber
&& + (-1)^{|z|(|x|+|y|)}\beta^{2n}(\{z,x,y\} \\\nonumber
&& =  \beta^{2n}[(x\ast y)\ast \alpha(z) +(-1)^{|x|(|y|+|z|)}(y\ast z)\ast \alpha(x)\\\nonumber
&& + (-1)^{|z|(|x|+|y|)}(z\ast x)\ast \alpha(y) +  \{x,y,z\} + (-1)^{|x|(|y|+|z|)}\{y,z,x\}\\\nonumber
&& + (-1)^{|z|(|x|+|y|)}\{z,x,y\}] \\\nonumber
& = &\beta^{2n}(0)\,\,\mbox{(by $(SHLY5)$ for $A_\alpha$)}\\\nonumber
& = & 0\nonumber
\end{eqnarray*}
 and thus we get $(SHLY5)$  for $A_{\beta^n}$. Next,\\
$\{x\ast_\beta  y,(\beta^n\alpha)(z),(\beta^n\alpha)(u)\}_\beta = \{\beta^{3n}(x\ast y),\beta^{3n}(\alpha(z)),\beta^{3n}(\alpha(u))\}\\ =
\beta^{3n}(\{x\ast y,\alpha(z),\alpha(u)\}).$
Therefore 
\begin{eqnarray*}  
&&\circlearrowleft_{(x,y,z)} (-1)^{|x||z|}\{x\ast_\beta  y,(\beta^n\alpha)(z),(\beta^n\alpha)(u)\}_\beta\\\nonumber
&& =  \circlearrowleft_{(x,y,z)}(-1)^{|x||z|}[\beta^{3n}(\{x\ast y,\alpha(z),\alpha(u)\})]\\\nonumber
&& =  \beta^{3n}(\circlearrowleft_{(x,y,z)}(-1)^{|x||z|}(\{x\ast y,\alpha(z),\alpha(u)\}))\\\nonumber
&& =  \beta^{3n}(0)\,\, \mbox{(by $(HLY6)$ for $A_\alpha$)}\\\nonumber
&& =  0
\end{eqnarray*}
So that we get $(SHLY6)$ for $A_{\beta^n}.$ Further, using $(SHLY7)$ for  $A_\alpha$ and  condition $\alpha\beta = \beta\alpha,$  we compute\\
 \begin{eqnarray*}
&&\{(\beta^n\alpha)(x),(\beta^n\alpha)(y),u\ast_\beta v\}_\beta =  \beta^{3n}(\{\alpha(x),\alpha(y),u\ast v\})\\
&& =  \beta^{3n}(\{x,y,u\}\ast \alpha^2(v) + (-1)^{|u|(|x|+|y|)}\alpha^2(u)\ast\{x,y,v\})\\
&& =  \beta^n(\beta^{2n}(\{x,y,u\})\ast(\beta^{2n}\alpha^2)(v)) + (-1)^{|u|(|x|+|y|)}\beta^n((\beta^{2n}\alpha^2)(u) \ast \beta^{2n}(\{x,y,v\}))\\
&& =  \{x,y,u\}_\beta \ast_\beta(\beta^{2n}\alpha^2)(v)  + (-1)^{|u|(|x|+|y|)}(\beta^{2n}\alpha^2)(u) \ast_\beta\{x,y,v\}_\beta\\
&& = \{x,y,u\}_\beta \ast_\beta(\beta^n\alpha)^2(v)  + (-1)^{|u|(|x|+|y|)}(\beta^n\alpha)^2(u) \ast_\beta\{x,y,v\}_\beta.
  \end{eqnarray*}

Thus $(SHLY7)$ holds for $A_{\beta^n}.$ Using repeatedly the condition $\alpha\beta = \beta\alpha$ and the identity $(SHLY8)$
for $A_\alpha,$ the verification for $(SHLY8)$for $A_{\beta^n}$ is as follows.\\
\begin{eqnarray*}
&  & \{(\beta^n\alpha)^2(x),(\beta^n\alpha)^2(y),\{u,v,w\}_\beta\}_\beta\\%\nonumber
& = & \beta^{2n}(\{(\beta^{2n}\alpha^2)(x),(\beta^{2n}\alpha^2)(y),\beta^{2n}(\{u,v,w\})\})\\%\nonumber
& = & \beta^{4n}(\{\alpha^2(x),\alpha^2(y),\{u,v,w\}\})\\%\nonumber
& = & \beta^{4n}(\{\{x,y,u\},\alpha^2(v),\alpha^2(w)\}) \\%\nonumber
& + & \beta^{4n}((-1)^{|u|(|x||y|)}\{\alpha^2(u),\{x,y,v\},\alpha^2(w)\})\\%\nonumber
& + & \beta^{4n}((-1)^{(|x|+|y|)(|u|+|v|)}\{\alpha^2(u),\alpha^2(v),\{x,y,w\}\})\\%\nonumber
& = & \beta^{2n}(\{\beta^{2n}(\{x,y,u\}),(\beta^{2n}\alpha^2)(v),(\beta^{2n}\alpha^2)(w)\})\\%\nonumber
& + & (-1)^{|u|(|x||y|)}\beta^{2n}(\{(\beta^{2n}\alpha^2)(u),\beta^{2n}(\{x,y,v\}),(\beta^{2n}\alpha^2)(w)\})\\%\nonumber
& + & (-1)^{(|x|+|y|)(|u|+|v|)}\beta^{2n}(\{(\beta^{2n}\alpha^2)(u),(\beta^{2n}\alpha^2)(v),\beta^{2n}(\{x,y,w\})\})\\%\nonumber
& = & \{\{x,y,u\}_\beta,(\beta^n\alpha)^2(v),(\beta^n\alpha)^2(w)\}_\beta\\%\nonumber
& + & (-1)^{|u|(|x||y|)}\{(\beta^n\alpha)^2(u),\{x,y,v\}_\beta,(\beta^n\alpha)^2(w)\}_\beta\\%\nonumber
& + & (-1)^{(|x|+|y|)(|u|+|v|)}\{(\beta^n\alpha)^2(u),(\beta^n\alpha)^2(v),\{x,y,w\}_\beta\}_\beta\\%\nonumber
\end{eqnarray*}

Thus $(SHLY8)$ holds for $A_{\beta^n}.$ Therefore, we get that $A_\beta$ is a Hom-LY superalgebra. This finishes the proof.
\hfill$\square$ \newline

In \cite{YAU2}, D. Yau etablished a general method of construction of Hom-algebras from their corresponding untwisted algebras.
From Theorem 3.1 we have the following method construction of Hom-LY superalgebras from LY superalgebras(this yields examples of Hom-Ly 
superalgebras). This method is an extension to binary-ternary superalgebras of  D. Yau's result (\cite{YAU2}, Theorem 2.3).
Such an extension to binary-ternary algebras is first mentioned in \cite{NOU1}, Corollary 4.5. \newline

\begin{Corollary}\label{Metho}
 Let $(A,\ast,[,,])$ be a LY algebra and $\beta$ an endomorphism of $(A,\ast,[,,]).$ If define on $A$  a
binary operation $"\tilde{\ast}"$ and a ternary operation $\{,,\}$ by\\
$x \tilde{\ast} y := \beta(x \ast y)$, \\
$\{x,y,z\}:= \beta^2([x,y,z]),$\\
then $(A,\tilde{\ast},\{,,\},\beta)$ is a Hom-LY algebra. \newline
\end{Corollary}

{\bf Proof.} The proof follows if observe that Corollary \ref{Metho} is Theorem \ref{Thm SHLY} when $\alpha = Id$ and $n = 1.$  
\hfill $\square$ \newline

\begin{Proposition}\label{H−SLY-Lie}
 Let $(L,[,]_{\alpha},\alpha)$ be a (multiplicative) Hom-Lie superalgebra. Define on $(L,[,]_{\alpha},\alpha)$ a ternary superoperation by
 \begin{equation}\label{T-Equa}
 \{x,y,z\}_{\alpha} := [[x,y],\alpha(z)].
 \end{equation}
 Then, $(L,[,]_{\alpha},\{,,\}_{\alpha},\alpha)$ is a Hom-Lie-Yamaguti Superalgebras.\newline
\end{Proposition}

\textbf{Proof:} Straightforward calculations by verification of $(SHLY1 - SHLY8)$ identities of Definition \ref{SHLY-Def} \hfill $\square$
%%%%%%%%%%%%%%%%%%%%%%%%%%%%%%%%%%%%%%%%%%%%%%%%%%%%%%%%%%%%%%%%%%%%%%%%%%%%%%%%%%%%%%%%%%%%%%%%%%%%%%%%%%%%%%%%%%%%%%%%%%%%%%%%%%%%%%%%%%%
\subsection{Examples of Hom-Lie-Yamaguti superalgebras}
%%%%%%%%%%%%%%%%%%%%%%%%%%%%%%%%%%%%%%%%%%%%%%%%%%%%%%%%%%%%%%%%%%%%%%%%%%%%%%%%%%%%%%%%%%%%%%%%%%%%%%%%%%%%%%%%%%%%%%%%%%%%%%%%%%%%%%%%%%

In the following, we give some examples of Hom-Lie-Yamaguti superalgebras which are constructed firstly from the example 2.8 in \cite{MAK5}
and using Proposition \ref{H−SLY-Lie}. In the second from the example 3.2 given in \cite{Eldu} (see also in \cite{MAK4}) and using the
Theorem \ref{Thm SHLY} via Corollary \ref{Metho} . Then, we obtain some family of  Hom-Lie-Yamaguti algebras of dimension 5 and dimension 3 
respectively. Indeed, 
\begin{Example}\label{EX-SLie}
 Consider the family of Hom-Lie superalgebra $osp(1,2)_{\lambda} = (osp(1,2),[,]_{\alpha_{\lambda}},\alpha_{\lambda})$ given in the 
 example 2.8 in \cite{MAK5}. The Hom-Lie superalgebra bracket $[,]_{\alpha_{\lambda}}$ on the basis elements is given, for $\lambda \neq 0,$ 
 by: 
 \begin{eqnarray*}
&& [H,X]_{\alpha_{\lambda}} = 2\lambda^2X, \quad  [H,Y]_{\alpha_{\lambda}} = \frac{-2}{\lambda^2}Y, \quad  [X,Y]_{\alpha_{\lambda}} = H,
 [Y,G]_{\alpha_{\lambda}} = \frac{1}{\lambda}F, \quad  [X,F]_{\alpha_{\lambda}} = \lambda G,\\
&& [H,F]_{\alpha_{\lambda}} = -\frac{1}{\lambda}F, \quad [H,G]_{\alpha_{\lambda}} = \lambda G, \quad
[G,F]_{\alpha_{\lambda}} = H, \quad  [G,G]_{\alpha_{\lambda}} = -2\lambda^2X, \quad  [F,F]_{\alpha_{\lambda}} = \frac{2}{\lambda^2}Y,
\end{eqnarray*}
where $\alpha_{\lambda} : osp(1,2) \rightarrow osp(1,2)$ is a linear map defined by
$$\alpha_{\lambda}(X) = \lambda^2X, \,\,\, \alpha_{\lambda}(Y) = \frac{1}{\lambda^2}Y, \,\,\, \alpha_{\lambda}(H) = H,\,\,\,
\alpha_{\lambda}(F) = \frac{1}{\lambda}F, \,\,\,\alpha_{\lambda}(G) = \lambda G, $$ and 
 $osp(1,2) = V_0 \oplus V_1$ is a Lie superalgebra where $V_0$ is generate by:
%  \begin{array}{cc}
 $$H =
 \left(
\begin{array}{ccc}
1 & 0 & 0\\
0 & 0 & 0\\
0 & 0 &-1
\end{array}
\right ),\,\,\,
X =
\left(
\begin{array}{ccc}
0 & 0 & 1\\
0 & 0 & 0\\
0 & 0 & 0
\end{array}
\right ),\,\,\,
Y =
\left(
\begin{array}{ccc}
0 & 0 & 0\\
0 & 0 & 0\\
1 & 0 & 0
\end{array}
\right )
$$
and $V_1$ is generated by:
$$
F =
\left(
\begin{array}{ccc}
0 & 0 & 0\\
1 & 0 & 0\\
0 & 1 & 0
\end{array}
\right ),\,\,\,
G =
\left(
\begin{array}{ccc}
0 & 1 & 0\\
0 & 0 & -1\\
0 & 0 & 0
\end{array}
\right ).
$$
such that 

$[H,X] = 2X, \quad  [H,Y] = -2Y, \quad  [X,Y] = H, \quad [Y,G] = F, \quad  [X,F] = G, $%\newline

$[H,F] = -F, \quad [H,G] = G, \quad [G,F] = H, \quad  [G,G] = -2X, \quad  [F,F] = 2Y,$\newline

Now, we define on $osp(1,2)_{\lambda} = (osp(1,2),[,]_{\alpha_{\lambda}},\alpha_{\lambda})$ a ternary superoperation by (\ref{T-Equa}).
Then, 
$sLY(1,2)_{\lambda} = (osp(1,2),[,]_{\alpha_{\lambda}},\{,,\}_{\alpha_{\lambda}},\alpha_{\lambda})$ is a family of Hom-Lie-Yamaguti 
superalgebras where the supercommutator $[,]_{\alpha_{\lambda}}$ is defined as above in this example \ref{EX-SLie} and the ternary 
superoperation $\{,,\}_{\alpha_{\lambda}}$ (we give only the ones with non zero values in the left hand side and using the identity $(SHLY4)$
of Definition \ref{SHLY-Def}, one can deduce the others values wich are non zero) is defined by 
\begin{eqnarray*}
 2H & = & \{H,X,Y\}_{\alpha_{\lambda}} = \{H,Y,X\}_{\alpha_{\lambda}} = 2\{H,F,G\}_{\alpha_{\lambda}} = 2\{H,G,F\}_{\alpha_{\lambda}}
= 2\{X,F,F\}_{\alpha_{\lambda}}\\
& = & -2\{Y,G,G\}_{\alpha_{\lambda}}  = -\{F,F,X\}_{\alpha_{\lambda}}\\ 
% & (= & -\{X,H,Y\}_{\alpha_{\lambda}} = -\{Y,H,X\}_{\alpha_{\lambda}} = -2\{F,H,G\}_{\alpha_{\lambda}} = -2\{G,H,F\}_{\alpha_{\lambda}}  \\
% &&  = -2\{F,X,F\}_{\alpha_{\lambda}} = 2\{G,Y,G\}_{\alpha_{\lambda}})\\
 -2\lambda^4X & = & \frac{1}{2}\{H,X,H\}_{\alpha_{\lambda}} =  -\{X,Y,X\}_{\alpha_{\lambda}} = \{F,G,X\}_{\alpha_{\lambda}} 
 = \{H,G,G\}_{\alpha_{\lambda}} = \{X,F,G\}_{\alpha_{\lambda}}\\
%  &(= & - \frac{1}{2}\{X,H,H\}_{\alpha_{\lambda}} = \{Y,X,X\}_{\alpha_{\lambda}} = -\{G,F,X\}_{\alpha_{\lambda}}= \{G,H,G\}_{\alpha_{\lambda}}
%  = \{F,X,G\}_{\alpha_{\lambda}}),\\
\frac{4}{\lambda^4}Y& =& \{H,Y,H\}_{\alpha_{\lambda}} = 2\{Y,X,Y\}_{\alpha_{\lambda}} = 2\{F,G,Y\}_{\alpha_{\lambda}} 
= -2\{H,F,F\}_{\alpha_{\lambda}} = 2\{Y,G,F\}_{\alpha_{\lambda}}\\
& = &  -\{F,F,H\}_{\alpha_{\lambda}}\\
% & (= & -\{Y,H,H\}_{\alpha_{\lambda}} = -2\{X,Y,Y\}_{\alpha_{\lambda}} = -2\{G,F,Y\}_{\alpha_{\lambda}} = 2\{F,H,F\}_{\alpha_{\lambda}}
% = -2\{G,Y,F\}_{\alpha_{\lambda}}),\\
-\frac{2}{\lambda^2}F & = & 2\{H,F,H\}_{\alpha_{\lambda}} = \{H,Y,G\}_{\alpha_{\lambda}} = 2\{H,G,Y\}_{\alpha_{\lambda}} 
= 2\{X,Y,F\}_{\alpha_{\lambda}} = -\{F,F,G\}_{\alpha_{\lambda}}\\
& = & 2\{X,F,Y\}_{\alpha_{\lambda}} = -2\{Y,G,H\}_{\alpha_{\lambda}} =  -2\{F,G,F\}_{\alpha_{\lambda}}\\
% &(= & -2\{F,H,H\}_{\alpha_{\lambda}} = -\{Y,H,G\}_{\alpha_{\lambda}} = -2\{G,H,Y\}_{\alpha_{\lambda}} = -2\{Y,X,F\}_{\alpha_{\lambda}}\\
% && = -2\{F,X,Y\}_{\alpha_{\lambda}} = 2\{G,Y,H\}_{\alpha_{\lambda}} = 2\{G,F,F\}_{\alpha_{\lambda}})\\
-\lambda^2 G & = & \{H,G,H\}_{\alpha_{\lambda}} = -\frac{1}{2}\{H,X,F\}_{\alpha_{\lambda}} = -\{H,F,X\}_{\alpha_{\lambda}} 
= -\{X,Y,G\}_{\alpha_{\lambda}} = \{X,F,H\}_{\alpha_{\lambda}}\\
& = & \{Y,G,X\}_{\alpha_{\lambda}} = \{F,G,G\}_{\alpha_{\lambda}}\\ 
% && (= -\{G,H,H\}_{\alpha_{\lambda}} = \frac{1}{2}\{X,H,F\}_{\alpha_{\lambda}} = \{F,H,X\}_{\alpha_{\lambda}} = \{Y,X,G\}_{\alpha_{\lambda}}
% = -\{F,X,H\}_{\alpha_{\lambda}}\\
% && \hspace{0.2cm}= -\{G,Y,X\}_{\alpha_{\lambda}}  = \{G,F,G\}_{\alpha_{\lambda}})
\end{eqnarray*}
These Hom-Lie-Yamaguti superalgebras are not Lie-Yamaguti superalgebras for $\lambda \neq \pm1.$ Indeed, the left hand side of the identity
$(SLY3)$, for $\alpha = Id,$ leads to $$\circlearrowleft_{(H,Y,X)}(-1)^{|H||X|}\left([H,Y],X] + \{H,Y,X\}\right) 
= 2\left(\frac{1-\lambda^4}{\lambda^2}\right)H$$
and also the left hand side of the identity $(SLY4)$, for $\alpha = Id,$ leads to
 $$\circlearrowleft_{(H,Y,X)}(-1)^{|H||X|} \{[H,Y],X,G\} = 2\left(1-\lambda^2\right)G.$$ Then, they not vanish for $\lambda \neq \pm1.$\newline
\end{Example}

\begin{Example}
From the example $M^3(3,1)$ of a non-Lie Malcev superalgebra given in \cite{Eldu} (see also in \cite{MAK4}) and defining on it a ternary 
superoperation by (\ref{M-LY}), we obtain a Lie-Yamaguti superalgebra $(sLY(3,1),[,],\{,,\})$ of dimension 4.  
defined with respect to a basis $(e_1,e_2,e_3,e_4)$, where $(sLY(3,1))_0 = span(e_1,e_2,e_3)$ and $(sLY(3,1))_1 = span(e_4)$, by the following
multiplication table %\newline
\begin{eqnarray*}
&& [e_1,e_3] = -e_1\quad (= -[e_3,e_1]), \,\,\,[e_2,e_3] = 2e_2 \quad(= -[e_3,e_2]),\,\,\,  [e_3,e_4] = -e_4 \quad(= -[e_4,e_3]),\\
&& [e_4,e_4] = e_1 + e_2,\\
&& \{e_1,e_3,e_3\} = 2e_1 \quad(= -\{e_3,e_1,e_3\}), \,\,\, \{e_2,e_3,e_3\} = 8e_2 \quad(= -\{e_3,e_2,e_3\}),\\
&& \{e_3,e_4,e_3\} = 2e_4 \quad(= -\{e_4,e_3,e_3\}),\,\,\, \{e_3,e_4,e_4\} = e_1 - 2e_2 \quad(= -\{e_4,e_3,e_4\})\\
&& \{e_4,e_4,e_3\} = -e_1 - 4e_2.
\end{eqnarray*}
Consider the even superalgebra endomorphism $\alpha_1: sLY(3,1) \rightarrow sLY(3,1)$ with respect to the same basis define by
 $$ \alpha_1(e_1) = a^2e_1,\quad \alpha_1(e_2) = a^2e_2, \quad \alpha_1(e_3) = be_1 + ce_2 + e_3, \quad \alpha_1(e_4) = a^2e_4 \quad \forall\,\,\, a,b,c \in \mathbb{K}$$
For each such even superalgebra endomorphism $\alpha_1$ and using the Theorem \ref{Thm SHLY} via Corollary \ref{Metho}, there is a 
Hom-Lie-Yamaguti superalgebra $sLY(3,1)_{\alpha_1} = (sLY(3,1),[,]_{\alpha_1},\{,,\}_{\alpha_1})$ such that
\begin{eqnarray*}
&& [e_1,e_3]_{\alpha_1} = -a^2e_1\quad (= -[e_3,e_1]_{\alpha_1}), \,\,\,[e_2,e_3]_{\alpha_1} = 2a^2e_2 \quad(= -[e_3,e_2]_{\alpha_1}),\\ 
&&[e_3,e_4]_{\alpha_1} = -ae_4 \quad(= -[e_4,e_3]_{\alpha_1}), \quad [e_4,e_4]_{\alpha_1} = a^2(e_1 + e_2),\\
&& \{e_1,e_3,e_3\}_{\alpha_1} = 2a^4e_1 \quad(= -\{e_3,e_1,e_3\}_{\alpha_1}), \,\,\, \{e_2,e_3,e_3\}_{\alpha_1} = 8a^4e_2 \quad
(= -\{e_3,e_2,e_3\}_{\alpha_1}),\\
&& \{e_3,e_4,e_3\}_{\alpha_1} = -2a^5e_4 \quad(= -\{e_4,e_3,e_3\}_{\alpha_1}),\,\,\, \{e_3,e_4,e_4\}_{\alpha_1} = a^4(e_1 - 2e_2) \quad
(= -\{e_4,e_3,e_4\}_{\alpha_1})\\
&& \{e_4,e_4,e_3\}_{\alpha_1} = a^4((2a-1)e_1 +2(a+1)e_2).
\end{eqnarray*}
While for example $$\circlearrowleft_{(e_3,e_4,e_4)}(-1)^{|e_3||e_4|}([[e_3,e_4],e_4] + \{e_3,e_4,e_4\}) = a^4(2(a-1)e_1 + (2a+3)e_2),$$
then for $a \neq 0,$ these Hom-Lie-Yamaguti superalgebras are note Lie-Yamaguti superalgebras. 

In the same a way, when we consider an even superalgebra endomorphism $\alpha_2: sLY(3,1) \rightarrow sLY(3,1)$ with respect to the same 
basis define by 
 $ \alpha_2(e_3) = be_1 + ce_2 +\frac{1}{2} e_3, \quad \alpha_2(e_4) = de_4, \quad \forall\,\,\, a,b,c,d \in \mathbb{K}$ and using the 
 Theorem \ref{Thm SHLY} via Corollary \ref{Metho}, we obtain a twisting of $sLY(3,1)$ into a family of Hom-Lie superalgebras define by
 $[e_3,e_4] = -de_2\,\,\, (= -[e_4,e_3])$ which are also Hom-Lie-Yamaguti superalgebras in particulary case where the ternary 
 superoperation is zero. Observe that they are also  Lie superalgebras and consequently they are Lie-Yamaguti superalgebras where the ternary
 superoperation is zero. 
 \end{Example}
 %%%%%%%%%%%%%%%%%%%%%%%%%%%%%%%%%%%%%%%%%%%%%%%%%%%%%%%%%%%%%%%%%%%%%%%%%%%%%%%%%%%%%%%%%%%%%%%%%%%%%%%%%%%%%%%%%%%%%%%%%%%%%%%%%%%%%%%%%
\section{\scshape{$n^{th}-$derived binary-ternary Hom-superalgebras}}
%%%%%%%%%%%%%%%%%%%%%%%%%%%%%%%%%%%%%%%%%%%%%%%%%%%%%%%%%%%%%%%%%%%%%%%%%%%%%%%%%%%%%%%%%%%%%%%%%%%%%%%%%%%%%%%%%%%%%%%%%%%%%%%%%%%%%%%%%%
In this section, we extend the notion of $n^{th}-$derived (binary) Hom-superalgebra to the case of ternary and binary-ternary
Hom-superalgebras. In particulary, we introduce  $n^{th}-$derived of Hom-Lie-Yamaguti superalgebra and we shown that the category
 of Hom-Lie-Yamaguti superalgebras is closed under the process of taking $n^{th}-$derived Hom-superalgebras.\newline
 
\begin{Definition}\label{H-T-D}
Let $A := (A,\{,,\},\alpha)$ be a ternary Hom-superalgebra and $n\geq 0$ an integer. Define on $A$ the $n^{th}-$derived ternary operation
$\{,,\}^{(n)}$ by
\begin{eqnarray}
\{x,y,z\}^{(n)} := \alpha^{{2^{n+1}}-2}(\{x,y,z\}), \forall\,\,  x,y,z \in A. 
\end{eqnarray}

Then $A^{(n)} := (A,\{,,\}^{(n)}, \alpha^{2^{n}})$ will be called the {\it $n^{th}-$one derived  ternary Hom-superalgebra} of A.
Now denote $\{x,y,z\}^{(n)} = \alpha^{{2^{n+1}}-2}(\{x,y,z\}).$ Then we note that
$A^0 = (A,\{,,\},\alpha), A^1 = (A, \{x,y,z\}^{(1)} = \alpha^2\circ \{x,y,z\},\alpha^2),$ and  $A^{n+1} = (A^n)^1$.\newline
\end{Definition}

\begin{Definition}\label{H-B-T-D}
 Let $A := (A,\ast,\{,,\},\alpha)$ be a binary-ternary Hom-superalgebra and $n\geq 0$ an integer. Define on $A$ the $n^{th}-$derived 
binary operation and the $n^{th}-$derived ternary operation $[,,]^{(n)}$ by
\begin{eqnarray}
&& x\ast^{(n)} y := \alpha^{{2^n}-1}(x\ast y)\\
&& \{x,y,z\}^{(n)} := \alpha^{{2^{n+1}}-2}(\{x,y,z\}), \forall\,\,  x,y,z \in A. 
\end{eqnarray}

Then $A^{(n)} := (A, \ast^{(n)},\{,,\}^{(n)}, \alpha^{2^{n}})$ will be called the {\it $n^{th}-$one derived  (binary-ternary) 
Hom-superalgebra} of A. Denote $\ast^{(n)} = \alpha^{2^{n}-1}\circ \ast$ and $\{x,y,z\}^{(n)} = \alpha^{{2^{n+1}}-2}(\{x,y,z\}).$ 
Then we note that
$A^0 = (A,\ast,\{,,\},\alpha),\,\, A^1 = (A, \ast^{(1)} = \alpha \circ \ast,\,\,\, \{x,y,z\}^{(1)} = \alpha^2\circ \{x,y,z\},\alpha^2),$ and 
$A^{n+1} = (A^n)^1$.
\end{Definition}
 One observes that, from Definition \ref{H-B-T-D}, if set $\{x,y,z\} = 0, \forall\,\,x, y, z \in A,$ we  recover the $n^{th}-$derived (binary)
 Hom-superalgebra of Definition \ref{H−B−D} and if $x \ast y = 0 \,\,\forall\,\,x, y \in A,$ then we have an $n^{th}-$derived
 (ternary) Hom-superalgebra (see Definition \ref{H-T-D})\newline
 
 \begin{Proposition}\label{H-LTS-D}
 Let $A := (A,\{,\},\alpha)$ be a multiplicative Hom-Lie supertriple system. Then for each $n\geq 0,$ the $n^{th}-$derived Hom-superalgebra 
 $A^{(n)} = (A,\{,,\}^{(n)} = \alpha^{{2^{n+1}}-2}\circ \{,,\},\alpha^{2^{n}})$ of $A$ is a multiplicative Hom-Lie supertriple system. %\newline
\end{Proposition}

\textbf{Proof:} The proof of this Proposition \ref{H-LTS-D} can be constitued throughout the proof of Theorem \ref{SHLYD} below if the binary 
superoperation is zero. $\hfill \square$ \newline

In the following result, we shown that the category of Hom-Lie-Yamaguti superalgebras is closed under the process of taking $n^{th}-$derived
Hom-superalgebras.
\begin{Theorem}\label{SHLYD}
Let $A_\alpha := (A,[,],\{, , \},\alpha)$ be a Hom-Lie-Yamaguti. Then, for each $n\geq 0$ $n^{th}$ derived Hom-algebra
$$A^{(n)} := (A,[,]^{(n)} = \alpha^{2^{n}-1} \circ [,], \{, , \}^{(n)}) = \alpha^{2^{n+1}-2} \circ \{,,\}, \alpha^{2^{n}}) $$ is a 
Hom-Lie-Yamaguti algebra. %\newline
\end{Theorem}

\textbf{Proof:} The identies $(SHLY1)- (SHLY4)$ for $A^{(n)}$ are obvious. The checking of $SHLY5$ for $A^{(n)}$ is as follows.
Indeed, consider $\circlearrowleft_{(x,y,z)}((-1)^{|x||y|}[[x,y]^{(n)},(\alpha^{2^{n}})(z)]^{(n)} + (-1)^{|x||y|}\{x,y,z\}^{(n)}).$  Then
\begin{eqnarray*}
 &&\circlearrowleft_{(x,y,z)}((-1)^{|x||z|}[[x,y]^{(n)},(\alpha^{2^{n}})(z)]^{(n)} + (-1)^{|x||z|}\{x,y,z\}^{(n)})\\\nonumber
 && =  \circlearrowleft_{(x,y,z)}((-1)^{|x||z|}[\alpha^{2^{n}-1}(\alpha^{2^{n}-1}([x,y])),(\alpha^{2^{n}-1}(\alpha(z)))]) 
 +  \circlearrowleft_{(x,y,z)}(\alpha^{2^{{n}+1}-2}((-1)^{|x||z|}\{x,y,z\}))\\\nonumber
&& = \alpha^{2^{{n}+1}-2}(\circlearrowleft_{(x,y,z)}((-1)^{|x||z|}[[x,y],\alpha(z)]) +  \circlearrowleft_{(x,y,z)}((-1)^{|x||z|}\{y,z,x\})) \\\nonumber
&& = \alpha^{2^{{n}+1}-2}(\circlearrowleft_{(x,y,z)}((-1)^{|x||z|}[[x,y],\alpha(z)] + (-1)^{|x||z|}\{x,y,z\}))\\\nonumber
&& = \alpha^{2^{{n}+1}-2}(0)\,\,\mbox{(by $(SHLY5)$ for $A_\alpha$)}\\\nonumber
&& = 0\nonumber
\end{eqnarray*}
 and thus we get $(SHLY5)$  for $A^{(n)}$.
 
 Next,
 \begin{eqnarray*}
\{[x,y]^{(n)},(\alpha^{2^{n}})(z),(\alpha^{2^{n}})(u)\}^{(n)} &= &\{\alpha^{2^{n}-1}([x,y]),(\alpha^{2^{n}})(z),(\alpha^{2^{n}})(u)\}^{(n)}\\
& = &\alpha^{2^{n+1}-2}(\{\alpha^{2^{n}-1}([x,y]),(\alpha^{2^{n}})(z),(\alpha^{2^{n}})(u)\})\\
& = & \alpha^{2^{n+1}-2}\alpha^{2^{n}-1}(\{[x,y],\alpha(z),\alpha(u)\})\\
& = & \alpha^{3(2^{n}-1)}(\{[x,y],\alpha(z),\alpha(u)\}).
\end{eqnarray*}
Therefore 
\begin{eqnarray*}  
&&\circlearrowleft_{(x,y,z)}((-1)^{|x||z|}\{[x,y]^{(n)},(\alpha^{2^{n}})(z),(\alpha^{2^{n}})(u)\}^{(n)})\\
& = &\,\, \circlearrowleft_{(x,y,z)}[\alpha^{3(2^{n}-1)}((-1)^{|x||z|}\{[x,y],\alpha(z),\alpha(u)\})]\\\nonumber
& = & \alpha^{3(2^{n}-1)}(\circlearrowleft_{(x,y,z)}((-1)^{|x||z|}\{[x,y],\alpha(z),\alpha(u)\}))\\\nonumber
& = &\alpha^{3(2^{n}-1)}(0) \,\,\mbox{(by $(SHLY6)$ for $A_\alpha$)}\\\nonumber
& = & 0
\end{eqnarray*}
So that we get $(SHLY6)$ for $A^{(n)}.$ Further, using $(SHLY7)$ for  $A_\alpha$ and  condition of multiplicativity linearity of $\alpha$, we
compute
 \begin{eqnarray*}
&&\{(\alpha^{2^{n}})(x),(\alpha^{2^{n}})(y),[u,v]^{(n)}\}^{(n)}\\ 
& = & \alpha^{2^{n+1}-2}(\{(\alpha^{2^{n}})(x),(\alpha^{2^{n}})(y),\alpha^{2^{n}-1}([u,v])\})\\
& = & \alpha^{2^{n+1}-2}\alpha^{2^{n}-1}(\{\alpha(x),\alpha(y),[u,v]\})\\
& = & \alpha^{2^{n+1}-2}\alpha^{2^{n}-1}([\{x,y,u\},\alpha^2(v)] +(-1)^{|u|(|x|+|y|)}[\alpha^2(u),\{x,y,v\}]) \\
& = & \alpha^{2^{n}-1}([\{x,y,u\}^{(n)},(\alpha^{2^{n}})^2(v)] +(-1)^{|u|(|x|+|y|)}[(\alpha^{2^{n}})^2(u),\{x,y,v\}^{(n)}])\\
& = & [\{x,y,u\}^{(n)},(\alpha^{2^{n}})^2(v)]^{(n)} +(-1)^{|u|(|x|+|y|)}[(\alpha^{2^{n}})^2(u),\{x,y,v\}^{(n)}]^{(n)}.
 \end{eqnarray*}
Thus $(SHLY7)$ holds for $A^{(n)} .$ Using repeatedly the condition of multiplicativity and the identity $(SHLY8)$
for $A_\alpha,$ the verification for $(SHLY8)$ for $A^{(n)}$ is as follows.
\begin{eqnarray*}
&  & \{(\alpha^{2^{n}})^2(x),(\alpha^{2^{n}})^2(y),\{u,v,w\}^{(n)}\}^{(n)}\\\nonumber
& = & \alpha^{2^{n+1}-2}(\{(\alpha^{2^{n}})^2)(x),(\alpha^{2^{n}})^2)(y),\alpha^{2^{n+1}-2}(\{u,v,w\})\}) \\\nonumber
& = & \alpha^{2^{n+1}-2}(\{(\alpha^{2^{n+1}})(x),(\alpha^{2^{n+1}})(y),\alpha^{2^{n+1}-2}(\{u,v,w\})\}) \\\nonumber
& = & (\alpha^{2^{n+1}-2})^2(\{\alpha^2(x),\alpha^2(y),\{u,v,w\}\})\\\nonumber
& = &  (\alpha^{2^{n+1}-2})^2(\{\{x,y,u\},\alpha^2(v),\alpha^2(w)\}) \\\nonumber
& + & (\alpha^{2^{n+1}-2})^2((-1)^{|u|(|x||y|)}\{\alpha^2(u),\{x,y,v\},\alpha^2(w)\})\\\nonumber
& + & (\alpha^{2^{n+1}-2})^2((-1)^{(|x|+|y|)(|u|+|v|)}\{\alpha^2(u),\alpha^2(v),\{x,y,w\}\})\\\nonumber
& = & (\alpha^{2^{n+1}-2})(\{(\alpha^{2^{n+1}-2})(\{x,y,u\}),(\alpha^{2^{n+1}-2})(\alpha^2(v)),(\alpha^{2^{n+1}-2})(\alpha^2(w))\})\\\nonumber
& + & (\alpha^{2^{n+1}-2})((-1)^{|u|(|x||y|)}\{(\alpha^{2^{n+1}-2})(\alpha^2(u)),(\alpha^{2^{n+1}-2})(\{x,y,v\}),(\alpha^{2^{n+1}-2})(\alpha^2(w))\})\\\nonumber
& + & (\alpha^{2^{n+1}-2})((-1)^{(|x|+|y|)(|u|+|v|)}\{(\alpha^{2^{n+1}-2})(\alpha^2(u)),(\alpha^{2^{n+1}-2})(\alpha^2(v)),(\alpha^{2^{n+1}-2})(\{x,y,w\})\})\\\nonumber
& = & \{\{x,y,u\}^{(n)},(\alpha^{2^{n}})^2(v),(\alpha^{2^{n}})^2(w)\}^{(n)}\\\nonumber
& + & (-1)^{|u|(|x||y|)}\{(\alpha^{2^{n}})^2(u),\{x,y,v\}^{(n)},(\alpha^{2^{n}})^2(w)\}^{(n)}\\\nonumber
& + & (-1)^{(|x|+|y|)(|u|+|v|)}\{(\alpha^{2^{n}})^2(u),(\alpha^{2^{n}})^2(v),\{x,y,w\}^{(n)}\}^{(n)}\\\nonumber
\end{eqnarray*}

Thus $(SHLY8)$ holds for $A^{(n)}.$ Therefore, we get that $A^{(n)}$ is a Hom-LY superalgebra. This finishes the proof. \hfill$\square$ \newline

\begin{Remark}
 \begin{itemize} 
  \item [$(1)$] If $[x,y]^{(n)} = 0,$ for all $x,y \in L,$ then  $n^{th}-$derived of Hom-Lie-Yamaguti superalgebras becomes an 
  $n^{th}-$derived  of Hom-Lie supertriple system  and, subsequently,  $n^{th}-$derived of ternary Hom-Nambu superalgebra 
  (since, by Definition \ref{LTSup}, any Hom-Lie supertriple system is automatically a ternary Hom-Nambu superalgebra). 
\item [$(2)$] If $\{x,y,z\}^{(n)} = 0,$ for all $x,y,z\in L,$ then the $n^{th}-$derived of Hom-LY superalgebra becomes a $n^{th}-$derived of 
Hom-Lie superalgebra. \newline
\end{itemize}  
\end{Remark}


\begin{thebibliography}{99}
\bibitem{MAK5}, K. Abdaoui, F. Ammar and A. Makhlouf, {\it Hom-Lie Superalgebras and Hom-Lie admissible Superalgebras}, Journal of algebra,
Vol. 324, Issue 7, 1513-1528 (2010).
\bibitem{MAK4}, K. Abdaoui, F. Ammar and A. Makhlouf, {\it Hom-alternative, Hom-Malcev and Hom-Jordan Superalgebras}, arXiv: 1304.1579v1 [math.RA]  Apr 2013.
\bibitem{Eldu} H. Albuquerque and A. Elduque, {\it Classification of Mal'tsev Superalgebras of small dimensions,} Algebra and Logic, Vol. 
\textbf{35} (6) (1996) 512-554.
\bibitem{Santana} H. Albuquerque A. P. Santana {\it Akivis superalgebras and speciality,} arXiv: 0710.4535v1 [math.RA] 24 oct 2007.

%\bibitem{Faouzi} F. Ammar, A. Makhlouf, {\it Hom-Lie superalgebras and Hom-Lie admissible Superalgebras,} arXiv: 0906.1668v2 [math.RA] 5 Aug 2009
\bibitem{Ataguema1} H. Ataguema, A. Makhlouf and S.D. Silvestrov, {\it Generalization of n-ary Nambu algebras and beyond,} J. Math. Phys., \textbf{50}
(2009), 083501.
\bibitem{ATTAN} S. Attan, A. Nourou Issa, {\it Hom-Bol algebras}, Quasigroups and related Systems \textbf{21} (2013), 131-146.
 \bibitem{Gohr2} Y. Fregier,  A. Gohr, {\it Unital algebras of Hom-associative type and surjective or injective twistings,} J. of Gen. Lie Theo.
and Appl. Vol. \textbf{3}, (2009), No. 4, 285-295.

% \bibitem {Ayupov} Sh. A. Ayupov and B. A. Omirov., {\it On Leibniz algebras,} in: Algebras and Operator Theory, Proceedings of the colloquium in Tashkent, 
% Kluwer Academic Publishers, Dordrecht (1998), 1-13.
% \bibitem{Benito} P. Benito, A. Elduque and F. Mat\'in-Herce, {\it Irreducible Lie-Yamaguti algebras,} J. Pure Appl. Alg., \textbf{213} (2009), 795-808.
% \bibitem{Casas} J. M. Casas, T. Pirashvili, {\it Ten-term exact sequences of Leibniz homology,} J. Algebra \textbf{231} (2000), 258-264.
% \bibitem{Felipe} R. Felipe, N. L\'opez-Reyes and F. Ongay, {\it R-Matrices for Leibniz algebras,} communicaci\'on T\'ecnica n$^\circ$ I-02-27/11-11-2002 
% (MB/CI MAT). 
\bibitem{GAP} D. Gaparayi and A. N. Issa, {\it A twisted generalization of Lie-Yamaguti algebras}, Int. J. Algebra \textbf{6} (2012), no. 7, 339-352.
\bibitem{HAR1} J. T. Hartwig, D. Larsson and S. D. Silvestrov, {\it Deformations of Lie algebras using $\sigma-$derivations,} J. Algebras,
 \textbf{292} (2006), 314-361. 
 \bibitem{NOU1} A. N. Issa, {\it Hom-Akivis algebras}, Comment. Math. Univ. Carolin. \textbf{52} (4) (2011), 485-500. 
 \textbf{292} (2006), 314-361. 
\bibitem{ANI} A. N. Issa, {\it Supercommutator algebras of right (Hom)-alternative superalgebras}, arXiv:1710.02706v1[math.RA]
 %\bibitem{Jacob} N. Jacobson, {\it Lie and Jordan triple systems,} Amer. J. Math., 71 (1949), 149-170.
  % \bibitem{Hofmann1} K. H. Hofmann and K. Strambach, {\it Lie's fundamental theorems for local analytic loops,} Pacific J. Math., \textbf{123} (1986), 
% 301-327.
% \bibitem{NOU3} A. N. Issa, {\it Remarks on the Construction of Lie-Yamaguti algebras from Leibniz algebras,} Int. J. Algebra, \textbf{5} (14) (2011, 
% 667 - 677.
% \bibitem{NOU2} A. N. Issa, {\it Some characterizations of Hom-Leibniz algebras,} arXiv:1011.1731v1.
%\bibitem{Kamiya} N. Kamiya, {\it On Lie algebras and triple systems},
%\bibitem{KIK1} M. Kikkawa, {\it Geometry of homogeneous Lie loops,} Hiroshima Math. J. \textbf{5} (1975), 141-179.
%\bibitem{KIN1} M. K. Kinyon and A. Weinstein, {\it Leibniz algebras, Courant algebroids and multiplications on reductive homogeneous spaces,} 
%  Amer. J. Math. \textbf{123} (2001), 525-550.
%  \bibitem{Loday1} J-L. Loday. {\it" Cyclic Homology"}, Grundlehren der Mathematischen Wissenschaften [Fundamental Principles of Mathematical Sciences],
%  301, Springer-Verlag, Berlin, 1992.
% \bibitem{Loday} J-L. Loday. {\it Une version non commutative des alg\`ebres de Lie: les alg\`ebres de Leibniz,} Enseign. Math., \textbf{39} (1993), 
% 269-293.
\bibitem{Frac} A. Koulibaly, M.F. Ouedraogo {\it Supersyst\`emes Triples de Lie associ\'es aux superalg\`ebres de Malcev} Africa Matematika, 
 s\'erie 3, Volume \textbf{14} (2002).
\bibitem{HAR2} D. Larsson and S. D. Silvestrov, {\it Quasi-Hom-Lie algebras, central extensions and 2-cycle-like identies,} J. Algebra
 \textbf{288} (2005), 321-344.
 \bibitem{HAR3} D. Larsson and S. D. Silvestrov, {\it Quasi-Lie algebras,} Comptemp.Math. \textbf{391} (2005).
% \bibitem{Lister} W. G. Lister, {\it A structure theory of Lie triple systems,} Trans. Amer. Math. Soc. SP (1952), 217 − 242.
\bibitem{MAK1} A. Makhlouf, {\it Hom-Alternative algebras and Hom-Jordan algebras,} Int. Elect. J. Alg., \textbf{8} (2010), 177-190.
%\bibitem{MAK2} A. Makhlouf, {\it Paradigm of nonassociative Hom-algebras and Hom-superalgebras,} Proceedings of Jordan structures in algebra and Analysis
Meeting, 143-177, Editional Circulo Rojo, Almeria, 2010.
\bibitem{MAK3} A. Makhlouf, Silvestrov S.D., {\it Hom-algebra structures,} J. Gen. Lie Theory Appl. \textbf{2} (2008), 51-64.
%\bibitem{MAK6} A. Makhlouf and S.D. Silvestrov, {\it Hom-algebras structures,} J. Gen. Lie Theory Appl., 2 (2008), 51-64.
%\bibitem{Nom1} K. Nomizu, {\it Invariant affine connections on homogeneous spaces,} Amer. J. Math., \textbf{76} (1954), 33-65.
% \bibitem{Pirashvili} T. Pirashvili, {\it On Leibniz homology,} Ann. Inst. Fourier (Grenoble), \textbf{44} (2) (1994), 401-411.
% \bibitem{Rakhimov1} I. S. Rakhimov and K. A. M. Atan, {\it On contractions and invariants of Leibniz algebras,} Bull. Malays. Math. Soc. (2) \textbf{35}
% (2A) (2012), 557-565.
% \bibitem{Okubo1} S. Okubo and N. Kamiya, {\it Quasi-classical Lie-super Algebra and Lie-super Triple Systems,} University of Rochester Report
% UR-1462, 1996.
\bibitem{Okubo2} S. Okubo, {\it Jordan-Lie Super Algebra and Jordan-Lie Triple System,} Journal of algebra \textbf{198}, 388-411 (1997).
% \bibitem{LY2} K. Yamaguti, {\it On the Lie triple system and its generalizations,} J. Sci. Hiroshima Univ., Ser. A 21 (1957/1958),155-160.
% 
% \bibitem{LY1}  K. Yamaguti, {\it On the Lie triple system and its generalization,} J. Sci. Hiroshima Univ., Ser A \textbf{21} (1957/1958), 107-113.
%\bibitem{YAU1} D. Yau, {\it Enveloping algebra of Hom-Lie algebras,} J. Gen. Lie Theory Appl., \textbf{2} (2008), 203-207.
\bibitem{YAU2} D. Yau, {\it Hom-algebras and homology,} J. Lie Theory, \textbf{19} (2009), 409-421.
\bibitem{YAU3} D. Yau, {\it Hom-Novikov algebras,} J. Phys. \textbf{A 44} (2011), 085202.
\bibitem{YAU4} D. Yau, {\it Hom-Maltsev, Hom-alternative and Hom-Jordan algebras,} International Electronic Journal of algebras, \textbf{11} (2012), 
 177-217.
 \bibitem{YAU5} D. Yau , {\it On n-ary Hom-Nambu and Hom-Nambu-Lie algebras,} J. Geom. Phys. \textbf{62} (2012), 506-522.
%%%%%%%%%%%%%%%%%%%%%%%%%%%%%%%%%%%%%%%%%%%%%%%%%%%%%%%%%%%%%%%%%%%%%%%%%%%%%%%%%%%%%%%%%%%%%%%%%%%%%%%%%%%%%%%%%%%%%%%%%%%%%%%%%%%%%%%%%%%%%%%%%%
\end{thebibliography}
\end{document}